\documentclass[12pt]{amsart}
\usepackage{geometry}                
\usepackage{graphicx}
\usepackage{amssymb}
\usepackage{amsmath}
\usepackage{epstopdf}
\usepackage{color}

\usepackage{url}
\usepackage{wasysym}
\usepackage[OT2, T1]{fontenc}
\usepackage[russian, english]{babel}
\title{On Arnold's and Pushkin's puzzles}

\author{Boris Khesin}
\address{Department of Mathematics,    University of Toronto, Canada;  \newline
e-mail: \texttt{khesin@math.toronto.edu}}

\begin{document}
\maketitle
\setcounter{footnote}{0}

\bigskip


A beautiful and thought provoking discovery by  Vladimir Arnold of a possible answer to Pushkin's puzzle
{\let\thefootnote\relax\footnotetext{\small *Editor's note for {\it Vladimir Arnold: Collected Works}, vol.VII, Springer Nature, 2025, pp.497-500.}}
stirred quite a bit of controversy and all kinds of opinions, both from Pushkinists and poetry lovers. 
Although this remarkable short note by Arnold [1] speaks for itself, we would like to draw the reader's attention
to the following passages that might have been dimmed in the English translation, but are quite noticeable 
and sufficiently prominent in the Russian original. 

\smallskip

The French epigraph is translated to Russian as follows:

\foreignlanguage{russian}{Proniknutyi0 t{shch}eslaviem, on obladal eshche toi0 osobennoi0 gordostp1yu, kotoraya pobuzhdaet priznavatp1sya s odinakovym ravnodushiem kak v svoikh dobrykh, tak i durnykh postupkakh,} --- \foreignlanguage{russian}{sledstvie chuvstva prevoskhodstva, bytp1 mozhet mnimogo.}\footnote{Steeped in vanity, he had even more of a special kind of pride that makes one acknowledge one's good and evil actions alike with the same indifference, out of a sense of superiority, perhaps an imaginary one. (transl. by T.Vardomskaya, this and the next footnote, see [1])}

Note those \foreignlanguage{russian}{toi0 osobennoi0 gordostp1yu} (a special kind of pride) and, particularly,    \foreignlanguage{russian}{bytp1 mozhet} (perhaps) at the end of a long phrase, 
thus supplying an additional degree of uncertainty to the description of the epigraph's protagonist. Later, in his note 
Arnold, when describing Pushkin's hoaxes, again employs a long phrase  with many subordinate clauses and with ``perhaps'' appearing in the middle, clearly mimicking Pushkin's style. 
\smallskip

And, remarkably, Arnold does it again in the description of his own, Arnold's unique research style: 

\foreignlanguage{russian}{
Buduchi po professii  [...] matematikom, 
ya vynuzhden v svoei0 rabote postoyanno opiratp1sya ne na dokazatelp1stva, a na oshchushcheniya, dogadki i gipotezy, perekhodya ot odnogo fakta k drugomu pri pomoshchi togo osobennogo vida 
ozareniya, kotoryi0 zastavlyaet usmatrivatp1 obshchie cherty 
v yavleniyakh, bytp1 mozhet, kazhushchikhsya po\-storonnemu vovse ne svyazannymi mezhdu soboi0.}\footnote{Being  [...] a mathematician, in my work I must
constantly depend not on proofs, but on sensations, guesses and hypotheses, moving from one fact to another by
means of a special kind of insight that lets one see commonalities in things that perhaps an observer may think completely unrelated.}

Of course, the latter phrase is a bit of trolling of the reader and archness on Arnold's part, 
since, as a mathematician, he founded most of his output on entirely rigorous proofs. Note, however, how in this 
ornate phrase with many subordinate clauses Arnold used  those ``a special kind'' and ``perhaps'' with an undoubted parallelism to the Onegin's florid lines. This pseudo-hidden playful attempt to use a similar to Puskin's style of writing only confirms Yuly Ilyashenko's
remarkable insight [2] that Arnold can be regarded as a Pushkin in mathematics.
\bigskip

Inspired by Arnold's love for trickery, puzzles and trolling, we would like to share with the reader  one
more finding in ``Eugene Onegin'', which Arnold would undoubtedly approve and enjoy (to the best of our knowledge
of his character). 
This recent discovery belongs to several scholars including Art\"em Stolpner from St.Petersburg, I.~Engelgardt, and B.~Bulatov [3],
and  their observation seemed to escape the attention of many  Pushkinists including Juri Lotman and Vladimir Nabokov.

\smallskip

Stolpner et al brought to light that Pushkin had hidden in his novel a hint that Tatiana and Olga Larin might have had different fathers!

\smallskip

Here is their reasoning. Recall the famous lines of Chapter 2 describing Praskovia Larina, Tatiana's mother:
\medskip

XXX\newline
\foreignlanguage{russian}{Ona lyubila Richardsona\newline 
Ne potomu, chtoby prochla, \newline
Ne potomu, chtob Grandisona\newline
 Ona Lovlasu predpochla; \newline
  No v starinu knyazhna Alina, \newline
  Ee0 moskovskaya kuzina, \newline
  Tverdila chasto ei0 ob nikh. \newline
  V to vremya byl eshche zhenikh \newline
  Ee0 suprug, no po nevole;\newline
   Ona vzdykhala po drugom, \newline
   Kotoryi0 serdtsem i umom \newline
   Ei0 nravilsya gorazdo bole: \newline
   Sei0 Grandison byl slavnyi0 frant, \newline
   Igrok i gvardii serzhant. }\newline
   
   XXXI \newline
 \foreignlanguage{russian}{  Kak on, ona byla odeta \newline
   Vsegda po mode i k litsu; \newline
   No, ne sprosyasp1 ee0 soveta, \newline
   Devitsu povezli k ventsu. \newline
   I, chtob ee0 rasseyatp1 gore, \newline
   Razumnyi0 muzh uekhal vskore \newline
   V svoyu derevnyu, gde ona, \newline
   Bog znaet kem okruzhena, \newline
   Rvalasp1 i plakala snachala, \newline
   S suprugom chutp1 ne razvelasp1; \newline
   Potom khozyai0stvom zanyalasp1, \newline
   Privykla i dovolp1na stala.}\footnote{The translations below and in the next footnote are due to J.E.Falen [4]
   (in his translation Princess' name is Laura): \newline

   30\newline
It wasn't that she'd read him, really,\newline
Nor was it that she much preferred\newline
To Lovelace Grandison, but merely\newline
That long ago she'd often heard\newline
Her Moscow cousin, Princess Laura,\newline
Go on about their special aura.\newline
Her husband at the time was still\newline
Her fianc\'e --- against her will!\newline
For she, in spite of family feeling,\newline
Had someone else for whom she pined --- \newline
A man whose heart and soul and mind\newline
She found a great deal more appealing;\newline
This Grandison was fashion's pet,\newline
A gambler and a guards cadet.\newline
   
   31\newline
   About her clothes one couldn't fault her;\newline
Like him, she dressed as taste decreed.\newline
But then they led her to the altar\newline
And never asked if she agreed.\newline
The clever husband chose correctly\newline
To take his grieving bride directly\newline
To his estate, where first she cried\newline
(With God knows whom on every side),\newline
Then tossed about and seemed demented;\newline
And almost even left her spouse;\newline
But then she took to keeping house\newline
And settled down and grew contented.   } \newline

To summarize, Tatiana's mother in her youth  was close to Princess Alina and before getting married she had what we would call today a boyfriend, ``Grandison''. As we learn from the novel, this romantic bonding was interrupted  by her marriage to Larin, after which her lifestyle changed abruptly. She soon gave birth to Tatiana, and later Olga. 
\smallskip

Now, let us compare this against a stanza from Chapter 7, when Tatiana and her mom come to Moscow for the ``bride fair'' and arrive at princess Alina's place:
\smallskip

XLI\newline
\foreignlanguage{russian}{--- Knyazhna,} mon ange! --- ``Pachette!'' --- \foreignlanguage{russian}{Alina! ---\newline
``Kto b mog podumatp1? Kak davno!\newline
  Nadolgo lp1? Milaya! Kuzina! \newline
  Sadisp1} --- \foreignlanguage{russian}{ kak e1to mudreno! \newline
  Ei0-bogu, stsena iz romana...''}\newline
 ---  \foreignlanguage{russian}{ A e1to dochp1 moya, Tatp1yana.} --- \newline
  \foreignlanguage{russian}{ ``Akh, Tanya! podoi0di ko mne} --- \newline
  \foreignlanguage{russian}{ Kak budto brezhu ya vo sne... \newline
  Kuzina, pomnishp1 Grandisona?''}\footnote{
  41\newline
`Princesse, mon ange!' `Pachette!' `Oh, Laura!'\newline
`Who would have thought?' `How long it's been!'\newline
`I hope you'll stay?' `Dear cousin Laura!'\newline
`Sit down. . . . How strange! . . . I can't begin . . .\newline
I'd swear it's from some novel's pages!'\newline
`And here's my Tanya.' `Lord, it's ages!\newline
Oh, Tanya sweet, come over here --- \newline
I think I must be dreaming, dear. . . .\newline
Oh, cousin, do you still remember\newline
Your Grandison?'  } \newline

Note the last four lines: as Tatiana approaches Alina,  the latter  cannot escape a deja vu.
It is virtually impossible to take these lines in any other way than assuming  that princess Alina recognizes  
Grandison's features in Tatiana's appearance! 
\smallskip

Admittedly, this can hardly be taken as proof. However, after having been directed to these lines in such a combination, it is difficult to forget about this possibility, as suggested by Stolpner et al. Pushkin definitely managed to 
seed a doubt in the attentive reader's mind. We believe that this curious and  spot-on discovery
would  be approved by Arnold, who liked finding such kinds of hints dispersed in the literature, and which, in a sense, nicely supplements Arnold's note on the epigraph. 

\bigskip
~
\bigskip

[1] V.Arnold: On the Epigraph to ``Evgenii Onegin'', -- {\it Izvestiya AN. Seriya Literatury i Yazyka}, vol. 56:2 (1997), 63; translated by Tamara Vardomskaya, see   
\url{https://vardomskaya.com/2016/11/27/598}

\smallskip

[2] Yu.Ilyashenko: V.I.Arnold, as I have seen him, -- in {\it ARNOLD: Swimming against the tide}, Eds: B.Khesin and S.Tabachnikov, AMS Providence, RI, 2014, 141-146.
\smallskip

[3] A.Stolpner: Personal communication (2019); I.Engelgardt: \url{https://otvet.mail.ru/question/19518323}, B.Bulatov: \url{https://bulatbayz.livejournal.com/993.html}
\smallskip

[4] A.Pushkin: Eugene Onegin, A Novel in Verse, -- translated by James E. Falen, Oxford World's Classics, 1995.

\end{document}